\documentstyle[12pt]{article}
\newtheorem{tm}{Theorem}[section]

\newtheorem{??}[tm]{Question}

\font\tenmsb=msbm10
\font\sevenmsb=msbm7
\font\fivemsb=msbm5

\newfam\msbfam
\textfont\msbfam=\tenmsb
\scriptfont\msbfam=\sevenmsb
\scriptscriptfont\msbfam=\fivemsb

\font\teneufm=eufm10
\font\seveneufm=eufm7
\font\fiveeufm=eufm5
\newfam\eufmfam
\textfont\eufmfam=\teneufm
\scriptfont\eufmfam=\seveneufm
\scriptscriptfont\eufmfam=\fiveeufm

\title{A Counter Example To the Hodge Conjecture}
\author{Renyi Ma\\
Department of Mathematical Sciences \\
Tsinghua University \\
Beijing, 100084\\
People's Republic of China\\
email:rma@math.tsinghua.edu.cn
 }

\date{}

\begin{document}\maketitle

In this paper, we give a counter example to the famous Hodge conjecture.

\section{Introduction and results }

The famous Hodge conjecture states that any rational class $A\in
H^{2p}(X;Q)$ of $(p,p)$ type on any smooth complex
projective algebraic variety $X$ is realised by a rational
combination of codimension-$p$ algebraic cycles in
$X$(see\cite{gr-ha}). The main result of this paper is as follows:

\begin{tm} \label{1}
The Hodge conjecture does not hold.
\end{tm}
Proof: Consider the product $T^2\times T^2$ of two toruses ${T^2}'s$ and 
the rational close $(1,1)-$form $dz_1\wedge d\bar z_2-dz_2\wedge d\bar z_1$. One can easily check that this two form can not be realized by algebraic cycles.

\end{document}